\documentclass[journal]{IEEEtran}
\usepackage[cmex10]{amsmath}
\usepackage[pdftex]{graphicx}
\usepackage{pdfsync}
\usepackage{amssymb}
\usepackage{color}
\usepackage{amsthm}
\usepackage{graphicx}
\usepackage{empheq}
\usepackage{hyperref}
\usepackage{array}
\usepackage{cite}

\newtheorem{prop}{Proposition}

\newtheorem{remark}{Remark}
\newtheorem{lem}{Lemma}

\def\f{\frac}
\def\pa{\partial}

\def\e{\eqref}

\def\i1n{i=1,\cdots,n}
\def\j1n{j=1,\cdots,n}
\def\ij1n{i,j=1,\cdots,n}

\def\R{\mathbb R}

\newcommand{\be}{\begin{equation}}
\newcommand{\ee}{\end{equation}}
\newcommand{\beq}{\begin{equation*}}
\newcommand{\eeq}{\end{equation*}}

\newtheorem{thm}{Theorem}

\begin{document}

\title{Feedback stabilization for the mass balance equations of a food extrusion process 
}

\author{
Mamadou~DIAGNE\thanks{Department of Mechanical and Aerospace Engineering,
University of California, San Diego,
La Jolla, CA, 92093, USA.
 E-mail: {\tt mdiagne@eng.ucsd.edu}.}
\ Peipei~SHANG\thanks{Department of Mathematics, Tongji University, Shanghai 200092, China.
E-mail: {\tt peipeishang@hotmail.com}.}
\
\ Zhiqiang~WANG\thanks{School of Mathematical Sciences and Shanghai Key Laboratory for Contemporary Applied Mathematics,
Fudan University, Shanghai 200433, China. E-mail: {\tt wzq@fudan.edu.cn}. }
}


\maketitle

\begin{abstract} 
In this paper, we study the stabilization problem for a food extrusion process in the isothermal case. 
The model  expresses the mass conservation in the extruder chamber and  consists of a hyperbolic Partial Differential  
Equation (PDE) and a nonlinear Ordinary  Differential Equation (ODE) whose dynamics describes the evolution of 
a moving interface.  
By using a Lyapunov approach, we obtain the exponential stabilization for the closed-loop system under
natural feedback controls through indirect measurements.
\end{abstract}


\IEEEpeerreviewmaketitle

\section{Introduction}

\IEEEPARstart{S}{crew} extruders  have become very popular for their ability to manufacture 
food and plastics products with  desired shapes and properties. Due to the   strong  interaction  between the mass, 
the energy and the momentum balances occurring  
in those processes, the  design of efficient controllers still remains a hard task at the industrial level. So far, the control oriented 
model of extruders are issued from some black box model of limited operational validity.  Following the objectives 
of the control, these models describes the extruder's  temperature and   flow  rate  at the die output  or the pressure 
dynamics based on   single input and single output or  multiple input and multiple output  system.
Generally,  extrusion processes are controlled using \emph{PID} \cite{MANOJ95, MMCA07, JUSTIN97}  
or predictive controllers  \cite{STUART00, Diagne13} with   oversimplified 
or empirical models.  In \cite{MANOJ95}, the volumetric expansion of the extrudate correlated to  the die 
temperature and  pressure and the  specific mechanical energy is 
chosen as the key product quality to be controlled. The authors  study the performance
of the \emph{PI}  controller based on the  regulation of the  die pressure using
feed rate as a manipulative variable and show that the
response of an improperly tuned controller may be too
sluggish on one hand, or too oscillatory on the other hand.  First-order, second-order and Lead-lag  Laplace
transfer-function are exploited in  \cite{ROSANA90} to design a feedforward controller  for  a twin-screw 
food extrusion process to reduce the effect
of feed rate and feed moisture content variations on the
die pressure. The  experimental results showed that  the die pressure response varies with the operation
 conditions, so a single gain value would not be suitable for all operating conditions. 
 The complexity of the process suggests adaptive control methods  because of the frequent changes occurring
during the extruder operation in the feed composition \cite{ROSANA90}. \cite{Lu1993} 
uses second-order  transfer functions while emphasizing  the difficulty in implementing those  types of model-based 
controllers due to the strong influence of   all the manipulated inputs and  measurable
process variables. Therefore,  very intelligent
controllers  need to be constructed  for extrusion
cooking process based on a control algorithm
developed from process experience.  
We mention that  the   transport delays, the strong interactions and the non-linearities make  it difficult to
control such systems   with   \emph{PID} controllers.  Predictive controllers might offer better 
performances but are somewhat difficult to implement  \cite{Correa93}.

 In the present work, we  consider the   stabilization of
  the  die   pressure to the desired setpoint in a food extrusion process. The controller is constructed based 
  on a  bi-zone model derived from the computation of 
   the conservation of mass in the extruder under the assumption of constant temperature and viscosity. 
   A Geometric decomposition of the extruder length in 
 \emph{ Partially  Filled  Zone (PFZ) } and \emph{ Fully  Filled  Zone (FFZ) }  allows to describe the 
 process by a  transport equation and a pressure-gradient equation defined on complementary time varying spatial domains. 
 The domains are coupled by a moving
interface whose dynamics is governed by an \emph{ODE} representing the total mass balance in
an extruder. We emphasize that an accurate die  pressure regulation is critical to ensure the uniformity 
of the extruded melt and is strongly related to the quality  of the final product. 
We propose a  suitable feedback control laws together
with practical measurements as output such that  the solution of the  closed-loop system
converges to a desired steady-state or equilibrium asymptotically.

The stabilization problems for hyperbolic systems has been widely studied in the literature.
The first approach  relies  on careful analysis of the classical solutions along
 the characteristics.  We refer to Greenberg and Li \cite{LiJMG} in the
case of second-order system of conservation laws and  more general situations on $n$th-order systems by Li in \cite{LiSt1}.

Another approach based on Lyapunov techniques was introduced by Coron et al. in \cite{CAB1}.
This approach was  improved by in \cite{CAB2} where  a strict Lyapunov function in terms of Riemann invariants
was constructed and  its time derivative can be made  negative definite by choosing properly the boundary conditions.
The Lyapunov function is very useful to analyze nonlinear hyperbolic systems of conservation laws because of its rubustness,
see \cite{CoronBook,  Coron08, DiagneCoron, CAB2, CWang13, XuCZ02, XuCZ09} for a wide range of applications to various models.
Among which, we are interested in a physical model for the
extrusion process which occurs very often in polymer material and food production.

The main contribution of this paper is to establish the exponential stabilization
 for the extrusion model under natural feedbacks by a Lyapunov function approach motivated by \cite{CAB2, CWang13, XuCZ02}.
The difficulties come from three aspects:  1) The domains on which the conservation laws are
defined depend on the solution through the dynamical interface;
2) The nonlinear coupling of the dynamics of the interface and the filling ratio
 in the \emph{Partially Filled Zone (PFZ)}  is not standard;
3) The measurements  and the feedbacks are natural, however,
some of the measurements are indirect, i.e.,
the measurements are not a part of solution but given through a indirect relation of solution.

The organization of the paper is as follows.
In Section \ref{sec2}, we give the physical description of the model.
The main result on stabilization (Theorem \ref{thm1})
and its  proof are given in Section \ref{sec3}.
Numerical simulations are provided in Section \ref{sec4}.
Finally in Section \ref{sec5},  we give our conclusion and some perspectives in control of the extrusion process.


\section{Description of the model}  \label{sec2}

\begin{figure}[h]
\begin{center}
\includegraphics[scale=0.5]{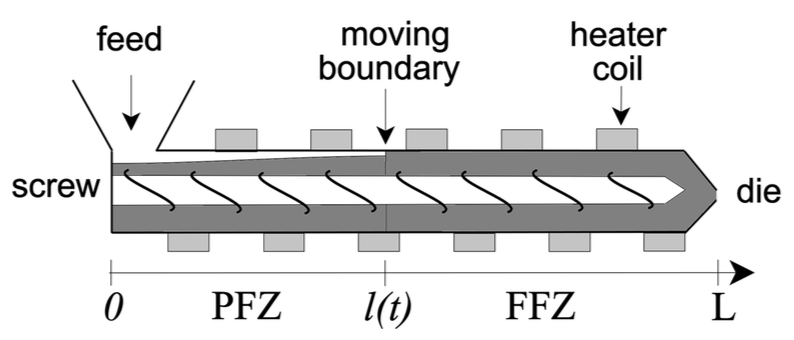}
\caption{A single-screw screw extruder} \label{fig1}
\end{center}
\end{figure}

An extruder is a process used for manufacturing objects
with fixed shapes and specific properties, see Fig. \ref{fig1}.  One or two Archimedean screws
are rotating inside the barrel in order to convect the extruded material from the feed to the die exit.
In food or polymers extrusion processes,  the ultimate control systems
involved manipulation of screw speed, feed rate, inlet solvent fraction, 
and barrel temperature for the regulation of moisture content, temperature and  viscosity of the finite product, residence time and die flow, etc.
In this paper, we consider the mass balances model \cite{Diagne, DSW1} motivated  by \cite{KULSH92,CHIN01} for cooking extrusion process.
In this case, the  material convection along  the  extruder chamber  of length $L$ is described in  two zones:  the  \emph{PFZ} ($[0,l(t)]$ in space) 
and a \emph{FFZ}  ($[l(t),L]$ in space) separated by a moving interface $l(t)$.  
The \emph{PFZ} and the \emph{FFZ} appear due to the die resistance that provokes an accumulation phenomena 
 and high pressure need to be  built-up to evict the extrudate out of the die.  
 By the mass conservation principle  the convection in the \emph{PFZ} is   
 described by the evolution of the filling ratio $f_p(t,x)$ for an homogeneous melt.
The melt  convection speed  in the \emph{PFZ}, namely, $\alpha_p$ depends 
on the screw speed $N(t)$ whereas the \emph{FFZ} transport velocity is related to the die pressure  
 $P(t,L)$: $\alpha_p \neq \alpha_f$.
Under the assumption of constant viscosity $\eta$ along the extruder,
the dynamics of the moving interface $l(t)$ is governed by  an \emph{ODE}
arising from the difference of the convection speed in the two regions. The flow rate in the \emph{FFZ} is constant and 
equal to the die flow rate $F_d(t)$  which is  proportional to the pressure difference
$\Delta P(t) := P(t,L)-P_0$ where  $P_0$ denotes  the atmospheric pressure.
For more detailed physical description of the model and definition of all the parameters appeared below,
one can refer to \cite{Diagne,DSW1}.

In this work, the  stabilization of  $(l(t),f_p(t,x))$ with the help of the actuated 
 screw speed $N(t)$ and inlet flow rate $F_{in}(t)$  is established 
based on  feedbacks that depend on the pressure difference $\Delta P(t)$ that is a practically useful measurement for the system. 
Considering the following  change of variables
\be \label{change-var}
x \mapsto y =\f{x}{l(t)}\  \text{in \emph{PFZ}} 
\  \text{and} \ 
x\mapsto y=\f{x-l(t)}{L-l(t)} \  \text{in \emph{FFZ}}, 
\ee
respectively, the time varying domains ($[0,l(t)]$, $[l(t),L]$) can be transformed to the fixed domain $[0,1]$ in space. 
For the sake of simplicity, we still denote by $x$ the space variable instead of $y$. 
More precisely, we consider the stabilization problem for the corresponding normalized
system on the spatial domain $[0,1]$. 
The interface dynamics which arises from a total mass balance  writes
\be\label{eqn-l}
\begin{cases}
\dot{l}(t)=\mathcal{F}(l(t),N(t),f_p(t,1)),\quad &\text{in}\ \mathbb{R}^+=(0,\infty),\\
l(0)=l^0,
\end{cases}
\ee
where
\be \label{F}
\mathcal{F}(l(t),N(t),f_p(t,1))=\frac{ \frac{K_d}{\eta}\Delta P(t)-\rho_{o}V_{eff}N(t)f_{p}(t,1)}{\rho_{o}S_{eff}(1-f_{p}(t,1))},
\ee
 $K_d$ , $\rho_0$, $\eta$ denote the die conductance, the melt density and the viscosity, respectively.  
 $V_{eff}$ and $S_{eff}$ are the effective volume and section between a screw element and the extruder  barrel, respectively.
\normalsize
Assuming a  constant viscosity along the extruder (the isothermal case),  the relation
\be \label{Delta P}
\Delta P(t)= \mathcal{P}(N(t),l(t))
:=\frac{\eta\rho_{o}
V_{eff}N(t) (L-l(t))}{B\rho_{o} +K_d
(L-l(t))}
\ee
is determined by integrating  the pressure-gradient equation corresponding to the momentum  balance in the \emph{FFZ} 
and considering a pressure continuity coupling relation at the normalized  interface, namely, $P(0,t)=P_0$ in the \emph{PFZ}  \cite{DSW1}.
The  filling ratio in the \emph{PFZ} writes
\be\label{eqn-fp}
\begin{cases}
\partial_t f_p(t,x)+\alpha_p \partial_x f_p(t,x)=0, \quad &\text{in}\ \mathbb{R}^+\times(0,1), \\
f_p(0,x)=f^{0}_p(x),\quad &\text{in}\ (0,1),\\
f_p(t,0)=\frac{F_{in}(t)}{\rho_{o} V_{eff}N(t)},\quad &\text{in}\ \mathbb{R}^+,
\end{cases}
\ee
where
\begin{align}\label{alpha-p}
\alpha_p=\f{\zeta N(t)-x \mathcal{F}(l(t),N(t),f_p(t,1))}{l(t)}.
\end{align}


\section{Main result and its proof } \label{sec3}

Let us define the constant equilibrium $l_e\in(0,L)$, $f_{p_e}\in(0,1)$, $N_e$,
$\Delta P_e$, $\alpha_{p_e}$ and $F_{in_e}$ by
$F(l_e,N_e,f_{p_e})=0$, $\Delta P_e=\mathcal{P}(l_e,N_e)$,
$ \alpha_{p_e}=\frac{\zeta N_e}{l_e}$, $F_{in_e}=\rho_{o}V_{eff}N_e f_{p_e}$.
Denote the difference
$\bar l(t):=l(t)-l_e, \bar N(t):=N(t)-N_e, \bar f_p(t,x):=f_p(t,x)-f_{p_e},  \bar F_{in}(t):=F_{in}(t)-F_{in_e},  
\Delta \bar P(t):=\Delta P(t)-\Delta P_e$ 
and the constants
\be \label{ai-bj}
\begin{cases}
 (a_1,a_2,a_3) = \Big ( \f{\pa F}{\pa l}, \f{\pa F}{\pa N}, \f{\pa F}{\pa f_p}\Big) \Big|_{(l_e,N_e,f_{p_e})},
\\
 (b_1,b_2)= \Big( \f{\pa \mathcal{P}}{\pa l}, \f{\pa \mathcal{P}}{\pa N}\Big )\Big|_{(l_e,N_e)}.
\end{cases}
\ee

The linear feedback law that we use is the following one:
\be\label{feed1}
\bar N(t)=k_1\cdot\Delta \bar P(t),\quad 
\bar F_{in}(t)=k_2\cdot\Delta \bar P(t),
\ee
where $\Delta P(t)$, thus $\Delta \bar P(t)$, is measurable.
The aim of stabilization is to find constants $k_1,k_2 \in \mathbb{R}$ such that  the closed-loop system
\eqref{eqn-l} and \eqref{eqn-fp} with feedback \eqref{feed1}
is asymptotically stable, i.e.,
$(\bar l(t), \bar f_p(t,\cdot)) \rightarrow 0 $ as $t\rightarrow\infty$.

Concerning the well-posedness of the Cauchy problem \eqref{eqn-l} and
\eqref{eqn-fp} with feedback \eqref{feed1}, we have the following proposition.
\begin{prop}\label{prop1}
Let $k_1,k_2\in \mathbb{R}$ be fixed.
There exists $\varepsilon>0$ such that for any $l^0\in \R$, $f^0_p\in H^2(0,1)$ satisfying
$ |l^0-l_e|^2+\|f^0_p(\cdot)-f_{p_e}\|^2_{H^2(0,1)}\leqslant \varepsilon$,
and the compatibility conditions at the point $(t,x)=(0,0)$,
system \eqref{eqn-l} and \eqref{eqn-fp}
with \eqref{feed1} has a unique solution
$(l,f_p) \in W^{1,\infty}([0,T)) \times C^0([0,T);H^2(0,1))$ for some $T\in (0,\infty]$. Moreover,
if
$ |l(t)-l_e|^2+\|f_p(t,\cdot)-f_{p_e}\|^2_{H^2(0,1)} \leqslant \varepsilon $ for all $ t\in [0,T),$
then $T=\infty$.
\end{prop}
\begin{remark}
The compatibility conditions at the point $(t,x)=(0,0)$ are the following:
\begin{align}\label{compa1}
 & f_p^0(0) = \f{F_{in}(0)}{\rho_{o} V_{eff} N(0)};
   \\
& \frac{\dot F_{in}(0)N(0)-F_{in}(0) \dot N(0) }{ \rho_{o} V_{eff} N(0)}
   +\f{\zeta N(0)}{l^0} f_p^{0'}(0)=0, \label{compa2}
\end{align}
where  $N(0), F_{in}(0)$ are determined by \eqref{Delta P}, \eqref{feed1} with $l(0)=l^0$,
while $\dot N(0)$,  $\dot F_{in}(0)$ are determined
 by differentiating  \eqref{Delta P} and  \eqref{feed1}  together with $\dot l(0)= F(l^0, N(0), f_p^0(1))$.
\end{remark}

The proof of Proposition \ref{prop1} is based on fixed point argument  and one can refer to
\cite{DSW1}  for the well-posedness of the corresponding open-loop system.
Our main result on stabilization of the interface  position $l(t)$ and the filling ratio $f_p(t,x)$ is the following theorem.
\begin{thm}\label{thm1}
Suppose that there exist  $k_1,k_2 \in \mathbb{R}$ such that
\begin{align}\label{k1}
&a_1+\f{k_1 a_2 b_1}{1-k_1 b_2}<0,  \\ \label{k2}
& \left |\f{a_3 b_1(k_2-f_{p_e}\rho_{o}V_{eff}k_1)}{\rho_{o}V_{eff}N_e(1-k_1b_2)} \right |
   <  \left |a_1+\f{k_1 a_2 b_1}{1-k_1 b_2} \right |,
\end{align}
where $a_1,a_2,a_3,b_1,b_2$ are given in \eqref{ai-bj}.
Then, the nonlinear system \eqref{eqn-l} and \eqref{eqn-fp}
 is locally exponentially stable under the feedback \eqref{feed1}, i.e.,
there exist constants $\varepsilon>0$, $M>0$ and $\omega>0$
such that for any $l^0\in \R$, $f^0_p\in H^2(0,1)$ satisfying
\be \label{small init}
|l^0-l_e|^2+\|f^0_p(\cdot)-f_{p_e}\|^2_{H^2(0,1)}\leqslant \varepsilon,
\ee
and the compatibility conditions at the point $(t,x)=(0,0)$,
the solution of \eqref{eqn-l} and \eqref{eqn-fp} with \eqref{feed1} satisfies
\begin{align}
& |l(t)-l_e|^2+\|f_p(t,\cdot)-f_{p_e}\|^2_{H^2(0,1)} \nonumber \\
& \leqslant  M e^{-\omega t}\Big(|l^0-l_e|^2+\|f^0_p(\cdot)-f_{p_e}\|^2_{H^2(0,1)}\Big),
\quad \forall t\geqslant 0.
\end{align}
\end{thm}

Before the proof of Theorem \ref{thm1}, we give several remarks.

\begin{remark}
The existence of $k_1,k_2\in \R$ satisfying
\eqref{k1}-\eqref{k2}  depends on $a_1,a_2,a_3, b_1,b_2$,
see Proposition  \ref{prop2} below. Its proof is given in Appendix.
\end{remark}

\begin{prop}\label{prop2}
There exist $k_1,k_2\in \R$ such that
\eqref{k1}-\eqref{k2}  hold
if and only if  either  1) $a_1<0$; or 2) $a_1\geqslant 0$ and $a_2b_1 \neq 0$.
\end{prop}

\begin{remark}
The measurement on $\Delta P(t)$ is of practical reason, thus the feedback \eqref{feed1} is indirect
in the sense that the measurements are made not on the solution $(l(t),f_p(t,x))$ itself.
\end{remark}

\begin{remark}
The proof of Theorem \ref{thm1} relies on a Lyapunov function approach.
The weight as $e^{-\gamma_i x}$ is essential to get a strict Lyapunov function. 
One can refer to the stabilization results by such weighted Lyapunov functions,
for quite general linear hyperbolic systems in \cite{XuCZ02,XuCZ09};
for one dimensional Euler equation with variable coefficients in \cite{GugatWang};
for a conservation law with nonlocal velocity in \cite{CWang13}.
\end{remark}

\noindent
{\bf Proof of Theorem \ref{thm1}:}
The construction of the Lyapunov functions is divided into three steps.

{\it Step 1. The stabilization of $l(t)$ and $f_p(t,\cdot)$ in $L^2(0,1)$.}

Let 
\begin{align}
V_0(t)=\bar l^2(t),  \quad 
\label{def-V1}
V_1(t)=\int_0^1 e^{-\gamma_1 x}\bar f^2_p(t,x)\, dx,
\end{align}
where $ \gamma_1>0$ is a constant to be  chosen later.

\begin{lem} \label{lem1}
There exist positive constants $A_1,\gamma_1,\beta_0,\beta_1,\delta_1$
such that the following estimates hold for every solution to system \eqref{eqn-l} and \eqref{eqn-fp}
with \eqref{feed1}
\begin{align} \label{dot V0V1}
 \dot V_0(t) +A_1 \dot V_1(t)  \leqslant & -(\beta_0+o(1)) V_0(t)
-(\beta_1+o(1)) V_1(t) \nonumber 
\\   & - (\delta_1+o(1)) \bar f_p^2(t,1),
\end{align}
where $o(1)$ represents various terms which tend to $0$ when
$|(\bar l(t),\bar N(t),\bar f_p(t,1)) |\rightarrow 0.$
\end{lem}

\noindent
{\bf Proof of Lemma \ref{lem1}:}
By definition of the equilibrium $(l_e,N_e,f_{p_e})$ and
the constants $(a_1,a_2,a_3,b_1,b_2)$, it is easy to get by expansion that
\begin{align}\label{est-F}
&\mathcal{F}(l(t),N(t),f_p(t,1))
=(a_1+o(1))\, \bar l(t)+(a_2+o(1))\, \bar N(t)\nonumber\\
&\hspace{36mm} +(a_3+o(1))\, \bar f_p(t,1),
   \\ \label{est-Delta P}
&\Delta\bar P(t)
=(b_1+o(1))\, \bar l(t)+(b_2+o(1))\, \bar N(t).
\end{align}
Furthermore, it follows from \eqref{feed1} and \eqref{est-Delta P} that
\begin{align}\label{est-Delta P2}
\Delta \bar P=&
\Big(\f{b_1}{1-k_1b_2}+o(1)\Big)\bar l(t),\\
\label{est-N}
\bar N(t)=&\Big(\f{k_1 b_1}{1-k_1b_2}+o(1)\Big)\bar l(t),\\
\label{est-Fin}
\bar F_{in}(t)=&\Big(\f{k_2 b_1}{1-k_1b_2}+o(1)\Big)\bar l(t).
\end{align}

Differentiating $V_0(t)$ with respect to $t$ and using  \eqref{eqn-l}, \eqref{est-F} and
\eqref{est-N}, one easily gets that
\begin{align} \label{est-V0}
\nonumber
\dot V_0(t)=& 2\bar l(t)\cdot \mathcal{F}(l(t),N(t),f_p(t,1))
\\ =&  (2 \theta_0+o(1) ) \bar l^2(t)
 + (2a_3+o(1))\bar l(t)\bar f_p(t,1).
\end{align}
where
$\theta_0 = a_1+\displaystyle \f{k_1 a_2 b_1}{1-k_1b_2} <0.$

On the other hand,
 \eqref{alpha-p} and \eqref{est-F} yield  that
\begin{align}\label{est-alpha px}
\alpha_p=\alpha_{p_e}+o(1), \quad
\alpha_{p_x}=o(1).
\end{align}
Differentiating $V_1(t)$ gives,
from \eqref{def-V1}, \eqref{est-alpha px}, \eqref{est-alpha px} and integration by parts, that
\begin{align}\label{est-V1}
\dot V_1(t)
=& BT_1+\int_0^1 (-\gamma_1\alpha_p  +\alpha_{p_x} )e^{-\gamma_1 x}\bar f^2_p(t,x)\, dx\nonumber\\
=& BT_1- (\gamma_1\alpha_{p_e}+o(1))V_1(t),
\end{align}
where
\be\label{BT1}
BT_1=-e^{-\gamma_1}(\alpha_{p_e}+o(1))\bar f^2_p(t,1)+(\alpha_{p_e}+o(1)) \bar f^2_p(t,0).
\ee
 Note that by \eqref{eqn-fp},\eqref{est-N}-\eqref{est-Fin}, we have
\begin{align}\label{est-fp0-2}
\bar f_p(t,0)
=\f{F_{in}(t)}{\rho_{o} V_{eff}N(t)}-f_{p_e} 
= (\theta_1+o(1) )\bar l(t),
\end{align}
where
$\theta_1= \f{b_1(k_2-f_{p_e}\rho_{o}V_{eff}k_1)}{\rho_{o}V_{eff}N_e(1-k_1b_2)}$.
Combining \eqref{est-alpha px}, \eqref{est-V1}, \eqref{BT1}, \eqref{est-fp0-2}, we get consequently
\begin{align}\label{est-V1-2}
\dot V_1(t)=& -(\gamma_1\alpha_{p_e}+o(1) )V_1(t)
-( e^{-\gamma_1} \alpha_{p_e} +o(1) )\bar f^2_p(t,1)\nonumber \\ 
& + (\alpha_{p_e} \theta_1^2 + o(1)) V_0(t).
\end{align}

Under the assumption of \eqref{k1}-\eqref{k2},
it is easy to get the existence of $A_1>0$ and $\gamma_1>0$  (suitably small) such that 
\be
   \left ( \begin{matrix} 
      2\theta_0 +A_1 \alpha_{p_e} \theta_1^2  &  a_3 \\
      a_3 & -A_1 e^{-\gamma_1} \alpha_{p_e} \\
   \end{matrix}
   \right)
\ee
is negative definite.
This  concludes the proof of   Lemma \ref{lem1} with \eqref{est-V0} and \eqref{est-V1-2}. \qed


{\it Step 2. The stabilization of $f_{p_x}(t,\cdot)$ in $L^2(0,1)$.}

Differentiating  \eqref{eqn-fp} with respect to $x$, we get
\be\label{eqn-fpx}
\begin{cases}
f_{p_{xt}}+\alpha_p f_{p_{xx}}+\alpha_{p_x}f_{p_x}=0,\\
f_{p_x}(0,x)=f^{0'}_p(x),\\
f_{p_x}(t,0)=\displaystyle\f{- f_{p_t}(t,0)}{\alpha_p |_{x=0}  },
\end{cases}
\ee
where
\be\label{def-fpt}
f_{p_t}(t,0)
=\f{\dot F_{in}(t)N(t)-F_{in}(t)\dot N(t)}{\rho_{o}V_{eff}N^2(t)}.
\ee
Let 
\be\label{def-V2}
V_2(t)=\int_0^1 e^{-\gamma_2x}f^2_{p_x}(t,x)\, dx,
\ee
where $ \gamma_2>0$ is a constant.

\begin{lem} \label{lem2}
There exist positive constants $\gamma_2,\beta_2, \delta_2,\theta_2$
such that the following estimate holds for every solution to system \eqref{eqn-l} and \eqref{eqn-fp}
with \eqref{feed1}
\begin{align} \label{dot V2}
\dot V_2(t) \leqslant & -(\beta_2+o(1)) V_2(t) - (\delta_2+o(1)) f_{p_x}^2(t,1)
 \nonumber\\ 
 &+\theta_2 (V_0(t) +\bar f_p^2(t,1) ).
\end{align}
\end{lem}

\noindent
{\bf Proof of Lemma \ref{lem2}.}
Differentiating \eqref{Delta P} and \eqref{feed1}  with respect to $t$ gives that,
\be\label{dot-feed}
\begin{cases}
\dot N(t)=k_1\cdot \Delta \dot P(t),\quad 
\dot F_{in}(t)=k_2\cdot \Delta \dot P(t),\\
\small  \Delta \dot P(t)=\f{\pa \mathcal{P}}{\pa l}\big|_{(l(t),N(t))} \dot l(t)
+ \f{\pa \mathcal{P}}{\pa N}\big|_{(l(t),N(t))} \dot N(t). 
\end{cases} 
\ee
Then it follows from \eqref{eqn-l}, \eqref{ai-bj}, \eqref{est-F}, \eqref{est-Delta P2} and \eqref{dot-feed} that
\begin{align}\label{dot-Delta P}
\Delta \dot P(t)=
&\f{\f{\pa \mathcal{P}}{\pa l}\big|_{(l(t),N(t))} }{1-k_1\f{\pa \mathcal{P}}{\pa N}\big|_{(l(t),N(t))}}\cdot \dot l(t)\nonumber\\
=&O(1)(\bar l(t)+\bar N(t)+\bar f_p(t,1)),
\end{align}
where $O(1)$ denotes various terms which are uniformly bounded
when $|(\bar l(t),\bar N(t),\bar f_p(t,1))| \rightarrow 0$.

Combining  \eqref{est-N}, \eqref{eqn-fpx}-\eqref{def-fpt} and \e{dot-feed}-\eqref{dot-Delta P}, we get easily  that
\begin{align}\label{fpx0}
f_{p_x}(t,0)
=&O(1)(\bar l(t)+\bar f_p(t,1)).
\end{align}
Differentiating \eqref{def-V2} results in,
by \eqref{est-alpha px}, \eqref{est-alpha px} and  \eqref{eqn-fpx}, that
\begin{align}\label{est-dot V2}
\dot V_2(t)
=&BT_2 +\int_0^1(-\gamma_2\alpha_p
      -\alpha_{p_x})e^{-\gamma_2x}f^2_{p_x}(t,x)\, dx\nonumber\\
=& BT_2+(-\gamma_2\alpha_{p_e}+o(1))V_2(t),
\end{align}
where
\begin{align}\nonumber 
BT_2=& (-e^{-\gamma_2}\alpha_{p_e}+o(1)) f^2_{p_x}(t,1) 
 \\&+(\alpha_{p_e}+o(1))f^2_{p_x}(t,0). \label{BT2} 
\end{align}
Thanks to \eqref{est-alpha px} and \eqref{fpx0}, \eqref{BT2} can be rewritten as
\begin{align}  \nonumber 
BT_2= &(-e^{-\gamma_2}\alpha_{p_e}+o(1))f^2_{p_x}(t,1)
 \\ &+O(1)(\bar l^2(t)+\bar f^2_p(t,1)),
\end{align}
which ends the proof of Lemma \ref{lem2} with  \eqref{est-dot V2}.
\qed

{\it Step 3. The stabilization of $f_{p_{xx}}(t,\cdot)$ in $L^2(0,1)$.}

By differentiating \eqref{eqn-fpx}, note that $\alpha_{p_{xx}}=0$, we derive that
\be\label{eqn-fpxx}
\begin{cases}
f_{p_{xxt}}+\alpha_p f_{p_{xxx}}+2\alpha_{p_x}f_{p_{xx}}=0,\\
f_{p_{xx}}(0,x)=f^{0''}_p(x),\\
f_{p_{xx}}(t,0)=\displaystyle\f{- f_{p_{xt}}(t,0)+\alpha_{p_x} |_{x=0} f_{p_x}(t,0)}{\alpha_p|_{x=0} },
\end{cases}
\ee
where
\begin{align}\label{fpxt0}
f_{p_{xt}}(t,0)=
\f{d}{dt}\Big(\f{-1}{\alpha_p |_{x=0} }\cdot \f{\dot F_{in}(t)N(t)-F_{in}(t)\dot N(t)}{\rho_{o}V_{eff}N^2(t)}\Big).
\end{align}

Let
\be\label{def-V3}
V_3(t):=\int_0^1 e^{-\gamma_3x}f^2_{p_{xx}}(t,x)\, dx,
\ee
where $\gamma_3>0$ is a positive constant.

\begin{lem} \label{lem3}
There exist positive constants $\gamma_3,\beta_3, \delta_3,\theta_3$
such that the following estimate holds for every solution to system \eqref{eqn-l} and \eqref{eqn-fp}
with \eqref{feed1}
\begin{align} \label{dot V3}
\dot V_3(t) \leqslant & -(\beta_3+o(1)) V_3(t) - (\delta_3+o(1)) f_{p_{xx}}^2(t,1)
 \nonumber\\ &+\theta_3 (V_0(t) +\bar f_p^2(t,1)+ f^2_{p_x}(t,1)).
\end{align}
\end{lem}

\noindent
{\bf Proof of  Lemma \ref{lem3}.}
Differentiating \eqref{def-V3} gives, by \eqref{est-alpha px}, \eqref{est-alpha px} and  \eqref{eqn-fpxx}, that
\begin{align}\label{est-dot V3}
\dot V_3(t)
=&BT_3 +\int_0^1\Big[-\gamma_3\alpha_p -3\alpha_{p_x}\Big]e^{-\gamma_3x}f^2_{p_{xx}}\, dx\nonumber\\
=&BT_3+  (-\gamma_3\alpha_{p_e}+o(1) )V_3(t),
\end{align}
where
\begin{align}\nonumber
BT_3=&-(e^{-\gamma_3}\alpha_{p_e}+o(1))f^2_{p_{xx}}(t,1) 
\\ \label{BT3} 
 &+(\alpha_{p_e}+o(1))f^2_{p_{xx}}(t,0).
\end{align}

In order to estimate $f_{p_{xx}}(t,0)$ or $f_{p_{xt}}(t,0)$,
essentially we need only to estimate $\ddot F_{in}(t)$ and $\ddot N(t)$,
according to  \eqref{dot-feed} and \eqref{fpxt0}.
On the other hand,
\eqref{dot-feed} simply yields that
\be\label{dot dot-feed}
\ddot N(t)=k_1\Delta \ddot P(t),\quad
\ddot F_{in}(t)=k_2\Delta \ddot P(t).
\ee
Therefore, from \eqref{eqn-l}, \eqref{est-F}, \eqref{dot-feed} and \eqref{dot-Delta P}, we have
\be \label{dot dot-Delta P}
\Delta \ddot P(t)
=O(1)\Big(\bar l(t)+\bar N(t)+\bar f_p(t,1)+f_{p_t}(t,1)\Big).
\ee
From \eqref{fpxt0} to \eqref{dot dot-Delta P}, we get
\begin{align}\label{est-fpxt0}
f_{p_{xt}}(t,0)
=O(1)\Big(\bar l(t)+\bar N(t)+\bar f_p(t,1)+f_{p_t}(t,1)\Big).
\end{align}
Combining \eqref{eqn-fp}, \eqref{est-N}, \eqref{fpx0}, \eqref{eqn-fpxx} and \eqref{est-fpxt0}, we get further
\begin{align}\label{est-fpxx0}
f_{p_{xx}}(t,0)
=&O(1)\Big(\bar l(t)+\bar f_p(t,1)+f_{p_x}(t,1)\Big).
\end{align}

By \eqref{est-alpha px} and \eqref{est-fpxx0}, \eqref{BT3} becomes
 \begin{align}  \nonumber 
BT_3= & -(e^{-\gamma_3}\alpha_{p_e}+o(1) )f^2_{p_{xx}}(t,1)
\\
&+O(1) (\bar l^2(t)+\bar f^2_p(t,1)+f^2_{p_x}(t,1)).
\end{align}
\normalsize
This implies the conclusion of Lemma \ref{lem3}.
\qed

{\it Step 4. The stabilization of $l(t)$ and $f_p(t,\cdot)$ in $H^2(0,1)$.}

Finally, let the Lyapunov function be
\be \label{def-L}
L(t)=V_3(t) + A_3 (V_2(t)+A_2 (V_0(t)+A_1V_1(t))),
\ee
where $A_1>0$ is such that \eqref{dot V0V1} holds and  $A_2,A_3>0$ will be chosen later.
Obviously, $L(t)$ is equivalent to $\bar l^2(t) +\|\bar f_p(t,\cdot)\|^2_{H^2(0,1)}$.
Then, by  \e{dot V0V1}, \e{dot V2}, \e{dot V3} and \eqref{def-L},
one can  choose $A_2>0$ and $A_3>0$ successively large enough such  that 
\be\label{est-dot L}
\dot L(t) \leqslant - (\beta+o(1)) L(t)
\ee
for some constant $\beta>0$.
We assume in a priori that
\be \label{small data}
 |(\bar l(t),\bar N(t),\bar f_p(t,1))| \leqslant \varepsilon_0
\ee
for some small $\varepsilon_0 >0$ such that
 $ |o(1)| \leqslant \frac{\beta}{2}$ in \e{est-dot L}.   Then $ \dot L(t)\leqslant - \f{\beta}{2} L(t)$,
thus $ L(t)\leqslant L(0) e^{-\f{\beta }{2} \, t}$. 
Thanks to the  assumption  \e{small init}, 
\e{small data} can be satisfied for all $t\geq 0$ if $\varepsilon>0$ is small enough.
The proof of Theorem \ref{thm1} is thus complete. \qed
%
%

\section{ Simulations} \label{sec4}

Computing the   time integration of the semi-discretized transport equations by finite volume
 with ODE45 routine of MATLAB, the stability result  
 is achieved under the assumptions of Theorems \ref{thm1}. 
 Especially, the existence of $k_1,k_2$ is guaranteed by 
  $a_1<0$ according to  Proposition  \ref{prop2}.  
  The gain  $k_1$ is chosen  to satisfy (\ref{k1}) and  (\ref{k2})   
 and  $k_2$ is derived from the  compatibility conditions (\ref{compa1})  and (\ref{compa2}) 
 for which the inlet flow rate $F_{in}(0)$ and the screw speed $N(0)$ with their respective time  derivatives
  $\dot {F}_{in}(0)$ and $\dot{N}(0)$ are computed with the help of  the feedback  (\ref{feed1}) and the initial filling ratio $f^0_p(0)$.  
\begin{itemize}
\item \textbf{Initial conditions:}\\
$ l^0=1.5 \, m$,  $f^0_p(x)=0.6905+0.025(1-\cos(\pi x))+ 0.0117 \sin(\pi x)$
\item \textbf{Setpoint values:}\\
$ l_e=1.37 \,m$,  $f_{pe}=0.6  $  
\item \textbf{Gain values:}\\
$ k_1=0.01$,  $ k_2=0.0001$ 
 \item \textbf{Verification of assumptions \eqref{k1}-\eqref{k2}}\\
$a_1= -0.0119 <0$
  \end{itemize}
\begin{figure}[htbp]
\centering
\begin{minipage}[t]{0.40\textwidth}
\centering
\includegraphics[scale=0.35]{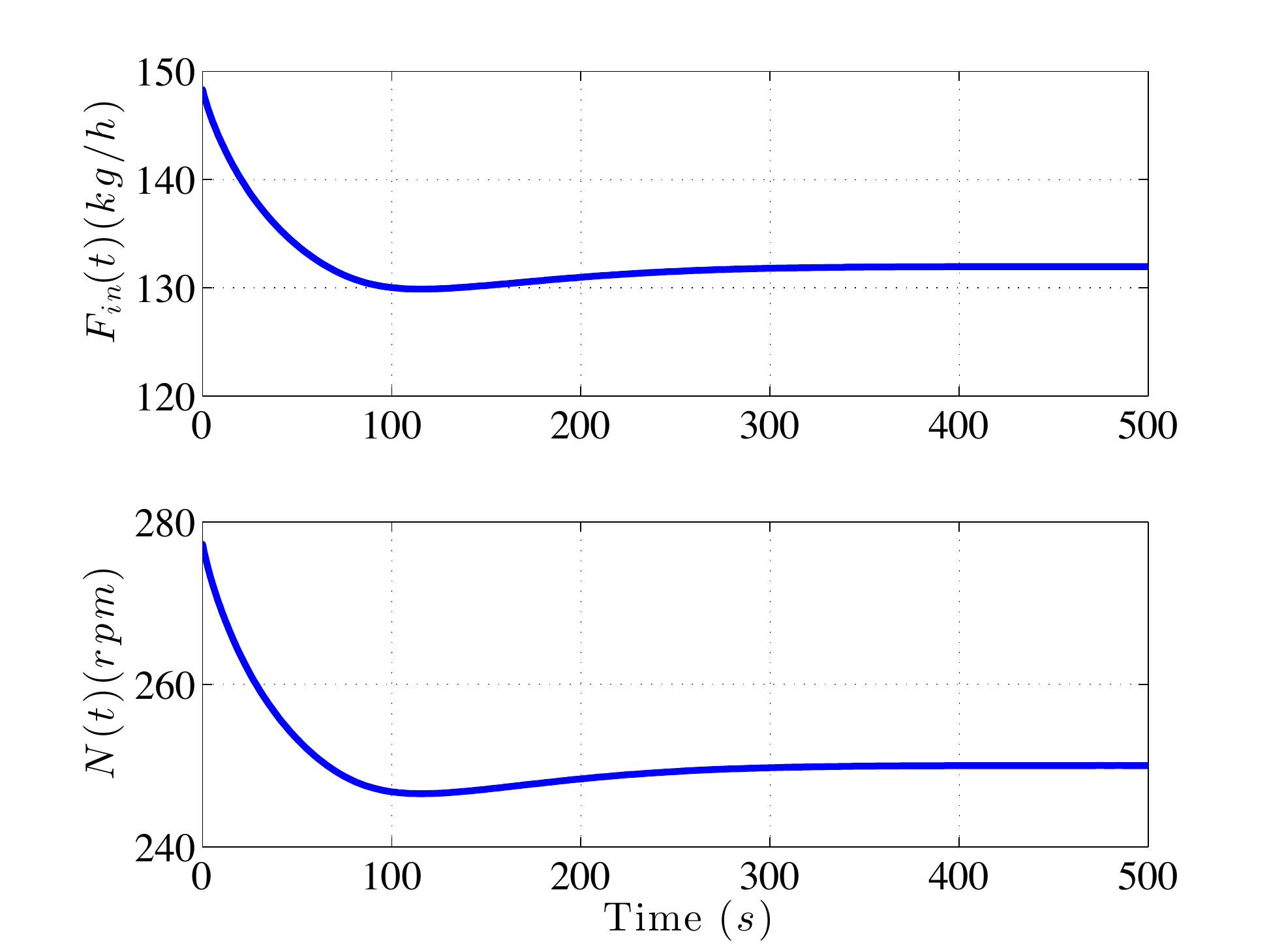}
\vspace{-5mm}
\caption{ Feed rate $F_{in}(t)$ and Screw speed $N(t) $ }
\end{minipage}
\begin{minipage}[t]{0.40\textwidth}
\centering
\includegraphics[scale=0.35]{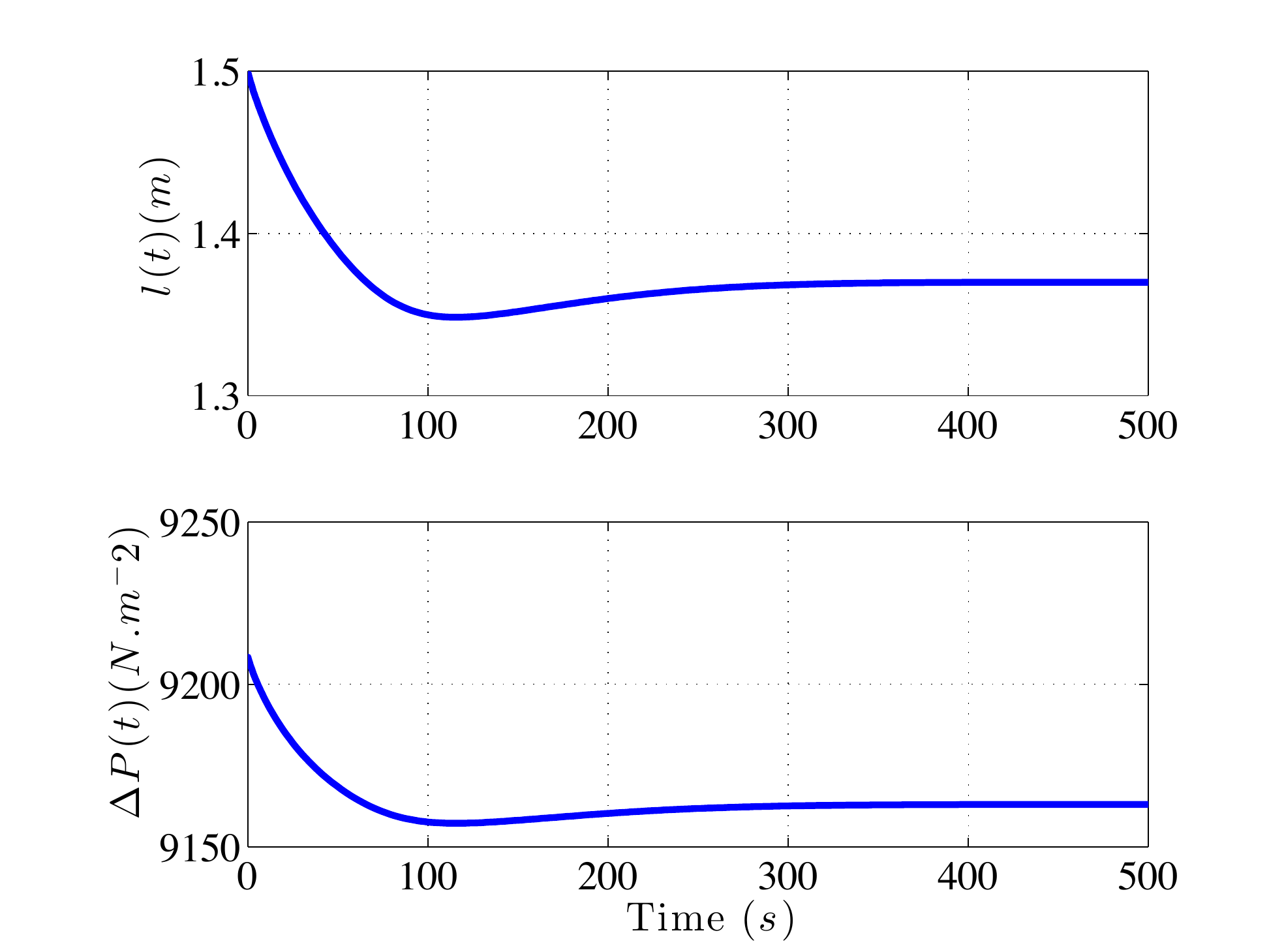}
\vspace{-5mm}
\caption{ Interface $l(t)$ and pressure difference $\Delta P(t)$ }
\end{minipage}
\begin{minipage}[t]{0.40\textwidth}
\centering
\includegraphics[scale=0.35]{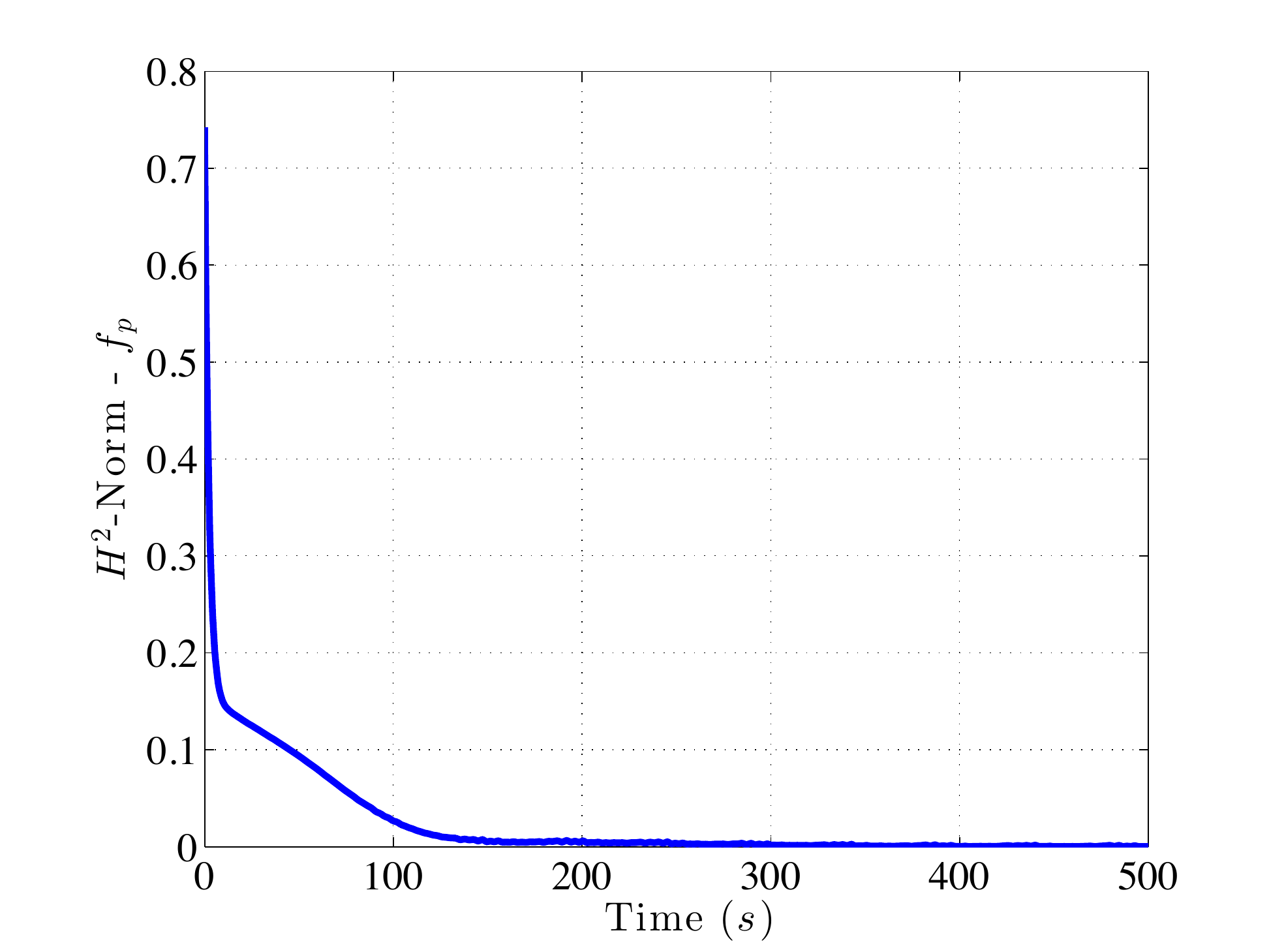}
\vspace{-5mm}
\caption{ Filling ratio $\|f_p(t,\cdot)-f_{p_{e}}\|_{H^2(0,1)}$}
\end{minipage}
\end{figure}

\section{Conclusion} \label{sec5}

In this paper, we study the stabilization  of
a physical model for the extrusion process,
which is described by  conservation laws  coupled through a dynamical interface.
The exponential stabilization is obtained for the closed-loop system with
natural feedback controls through indirect measurements.
The proof relies on Lyapunov approach.
Numerical simulations are made as supplementary to the theoretical results.
As a future work, it is interesting  to study the controllability of boundary profile, i.e., to reach the desired moisture
and temperature at the die under suitable controls.
This problem is rather challenging for mathematical theory but also  very useful in applications. 
 The proposed result is a first step towards  controlling complex screw  
 extrusion systems that include unknown parameters and time-varying  disturbances acting on the pressure dynamics and the viscosity in the FFZ. 


%

\appendix

\small
\noindent
{\bf Physical definition of the parameters}\\
\begin{tabular}{ll}
$L=2$  $m$     &\text{Extruder Length} \\
$B=2.4\times10^{-6}$    $m^{4}$           &\text{Geometric parameter}  \\
 $K_d=2.4\times10^{-3}$  $m^{3}$    &\text{Geometric parameter} \\
  $\zeta=0.003$  $m$     &\text{Screw Pitch }\\
$\eta=125$  $Pa$    $s^{-1}$  & \text{Melt viscosity} \\
 $\rho_{o}=350$  $kg$    $m^{-3}$  &\text{Melt density} \\
$S_{eff}=0.014$ $m^2$ &  \text{Effective area} \\
$V_{eff}=\zeta S_{eff}$  &  \text{Effective volume}\\
\end{tabular}
\normalsize
\vspace{5mm}

\noindent
{\bf Proof of Proposition \ref{prop2}.} It suffices to consider the existence of $k_1$ satisfying \eqref{k1}.
\it{Case 1) $a_1<0$. }  \eqref{k1} is true for $k_1 \rightarrow 0$.
\it{Case 2) $a_1\geqslant 0$} and  $a_2b_1\neq 0$.
If $b_2=0$, \eqref{k1} is true either for $k_1 \rightarrow +\infty$ or for $k_1 \rightarrow -\infty$;
If $b_2\neq 0$, \eqref{k1} is true either for  $k_1\rightarrow \displaystyle\f{1}{b_2} + $
or for $k_1\rightarrow \displaystyle \f{1}{b_2} -$. Upon the existence of $k_1$,
\eqref{k2} is true  either for $k_2$ is arbitrary if $a_3b_1=0$
or for $k_2 \rightarrow f_{p_e}\rho_{o}V_{eff}k_1$ if  $a_3b_1 \neq 0$.
\qed



\section*{Acknowledgment}
\small
The authors would like to thank Professor Jean-Michel Coron and
Professor Miroslav Krstic for their helpful comments and constant support.
The authors are thankful to  the support  of the ERC advanced grant 266907 (CPDENL)
and the hospitality of the Laboratoire Jacques-Louis Lions of Universit\'{e} Pierre et Marie Curie.
Peipei Shang was partially supported by the National Science Foundation of China (No. 11301387)
and by the Shanghai Pujiang Program (13PJ1408500). Mamadou Diagne  is  currently supported by
the Cymer Center for Control Systems and  Dynamics  of University of California San Diego as a postdoctoral fellow.
Zhiqiang Wang was partially supported by the National Science Foundation of China (No. 11271082)
and by the State Key Program of National Natural Science Foundation of China (No. 11331004).

\normalsize

\ifCLASSOPTIONcaptionsoff
  \newpage
\fi



%


\bibliographystyle{plain}
\bibliography{s}
\end{document}